\documentclass[12pt]{report}
\usepackage[brazil]{babel}
\usepackage[latin1]{inputenc}
\usepackage{hyperref}
\usepackage{amsmath,amssymb,indentfirst}

\pdfoutput=1

\begin{document}
\centerline{{\bf THE FUNDAMENTAL THEOREM OF ALGEBRA:}}
\centerline{{\bf A MOST ELEMENTARY PROOF }}

\

\centerline{{\bf Oswaldo Rio Branco de Oliveira}}

\

\hspace{- 0,5 cm}{\bf Abstract.} This paper shows an elementary and direct proof of the Fundamental Theorem of Algebra, via Bolzano-Weiestrass Theorem on Minima and the Binomial Formula, that avoids: any root extraction other than the one used to define the modulus function over $\mathbb C$, trigonometry, differentiation, integration, series, arguments by induction and $\epsilon\hspace{-0,1 cm}-\hspace{-0,1 cm}\delta$ type arguments.

\vspace{0,2 cm}

\hspace{- 0,5 cm}{\sl Mathematics Subject Classification: 12D05, 30A10}

\hspace{- 0,5 cm}{\sl Key words and phrases:} Fundamental Theorem of Algebra, Inequalities in $\mathbb C$. 

The aim of this article is, by combining an inequality proved in [9] and a lemma by Estermann [4], to show a very elementary proof of the FTA that requires no other nth root than the square root implicit in the modulus function. Following a suggestion given by Littlewood [8], see also Remmert [11], the proof requires a mininum amount of ``limit processes lying outside algebra proper''. Hence, the proof avoids differentiation, integration, series, angle and the transcendental functions (i.e., non-algebraic functions)
$\cos \theta$, $\sin\theta$ and $e^{i\theta}$, $\theta \in \mathbb R$. Another reason to avoid these functions is justified by the fact that the theory of transcendental functions is more profound than that of the FTA (a polynomial result), see Burckel [3]. Also avoided are arguments by induction and $\epsilon-\delta$ type arguments.

Many elementary proofs of the FTA, implicitly assuming the modulus function $|z|=\sqrt{z\overline{z}}$, where $z \in \mathbb C$, assume either the Bolzano-Weierstrass Theorem on Minima or the Intermediate Value Theorem, plus polynomial continuity. Then, along the proof it is used further root extraction in $\mathbb R$ or in $\mathbb C$ (see Argand [1] and [2], Estermann [4], Fefferman [5], Kochol [7], Littlewood [8], Oliveira [9], Redheffer [10], Remmert [11], Searcóid [12], Vaggione [13]). Beginning with Littlewood [8], some of these proofs include a proof by induction of the existence of every nth root, $n\in \mathbb N$, of every complex number (see [7], [11] and [12]). 

\newpage

Besides the modulus function (derived from the inner product $z\overline{w}$, with $z,w\in \mathbb C$) and the Binomial Formula $(z+w)^n= \sum_{j=0}^n\binom{n}{j}z^jw^{n-j}$, $z\in \mathbb C$, $w\in \mathbb C$, $n\in \mathbb N$, $\binom{n}{j}=\frac{n!}{j!(n-j)!}$ and $0!=1$, it is assumed, without proof, only: 
\begin{itemize}
\item[$\bullet$] Polynomial continuity.
\item[$\bullet$] Bolzano-Weierstrass Theorem: {\sl Any continuous function $f:D \to \mathbb R$, $D$ a bounded and  closed disc,
has a minimum on $D$}. 
\end{itemize}

Right below we show, for the case $k$ even, $k\geq 2$, a pair of inequalities that Estermann [4] proved for every $k\in \mathbb N \setminus\{0\}$. The proof, via binomial formula, is a simplification of the one made by induction and given by Estermann. The case $k$ odd can be proved similarly, if one wishes.

\vspace{0,2 cm}

\hspace{-0,6 cm}{\bf Lemma (Estermann).} For $\zeta = \Big(1+\frac{i}{k}\Big)^2$ and $k$ even, $k\geq 2$, we have
\[\textrm{Re}[\zeta^k]\ <\ 0 \ <\ \textrm{Im}[\zeta^k]\ .\] 
{\bf Proof.} Since $k=2m$ and $2k=4m$, for some $m \in \mathbb N$, applying the formulas
\begin{displaymath}
\begin{array}{ll}
\textrm{Re}\left[\left( 1 + \frac{i}{k}\right)^{2k}\right] &= 1 -\binom{2k}{2}\frac{1}{k^2} +\binom{2k}{4}\frac{1}{k^4} \,+\,\sum\limits_{\textrm{odd}\,j\,, j=3}^{k-1}\left[ -\binom{2k}{2j}\frac{1}{k^{2j}}+\binom{2k}{2j+2}\frac{1}{k^{2j+2}}\right]\ \textrm{and}\,,\\
\textrm{Im}\left[\left( 1 + \frac{i}{k}\right)^{2k}\right] & = \sum\limits_{\textrm{odd}\, j\,, j=1}^{k-1}\left[\binom{2k}{2j-1}\frac{1}{k^{2j-1}} \,-\,\binom{2k}{2j+1}\frac{1}{k^{2j+1}}\right]\ ,
\end{array}
\end{displaymath}
we end the proof by noticing that for every $j\in \mathbb N$, $1\leq j\leq k-1$, we have
\begin{displaymath}
\begin{array}{ll}
1 -\binom{2k}{2}\frac{1}{k^2} +\binom{2k}{4}\frac{1}{k^4} & =  1 - \Big(2 -\frac{1}{k}\Big)\Big(\frac{2}{3} + \frac{5}{6k} -\frac{1}{2k^2}\Big) \leq \\
&\leq 1 - \frac{3}{2}\Big(\frac{2}{3} +\frac{5k-3}{6k^2}\Big)=-\frac{3}{2}\cdot\frac{5k-3}{6k^2}<0 \,,\\
\\
-\binom{2k}{2j}\frac{1}{k^{2j}}+\binom{2k}{2j+2}\frac{1}{k^{2j+2}} & =-\frac{(2k)!}{(2j)!\,k^{2j}\,(2k-2j-2)!}\left[\frac{1}{(2k-2j)(2k-2j-1)} -\frac{1}{(2kj+2k)(2kj+k)}\right]<0\,,\,\\
\\
\binom{2k}{2j-1}\frac{1}{k^{2j-1}} -\binom{2k}{2j+1}\frac{1}{k^{2j+1}} & =\frac{(2k)!}{(2j-1)!\,(2k - 2j -1)!}\frac{1}{k^{2j-1}}\left[\frac{1}{(2k -2j + 1)(2k -2j)} - \frac{1}{(2kj +k)(2kj)}\right] > 0\ \ .
\end{array}
\end{displaymath}
\vspace{0,2 cm}

\newpage

\hspace{-0,6 cm}{\bf Theorem.} Let $P$ be a complex polynomial, with degree$(P)=n\geq 1$. Then, $P$ has a zero.

\hspace{-0,6 cm}{\bf Proof.} Putting $P(z)=a_0+a_1z+...+a_nz^n$, where $a_j\in \mathbb C$, $0\leq j\leq n$, $a_n \neq 0$, we  have $|P(z)|\geq |a_n||z|^n-|a_{n-1}||z|^{n-1}-...-|a_0| |z|^n|$. Hence,   $|P(z)|\rightarrow \infty$ as $|z|\rightarrow \infty$ and, by continuity, $|P|$ has a global minimum at some $z_0\in \mathbb C$. We can suppose without loss of generality $z_0=0$. Hence, 
\[(1)\ \ \ \ \ \ \ \ \ \ \ \ \ \ \ \ \ \ \ \ \ \ \ \ \ \ \ |P(z)|^2 - |P(0)|^2 \geq 0\,,\ \forall\, z \in \mathbb C\,,\ \ \ \ \ \ \ \ \ \ \ \ \ \ \ \ \ \ \ \ \ \  \ \ \ \ \ \ \ \ \ \ \ \ \ \ \ \ \ \]

\hspace{-0,6 cm}and $P(z) = P(0)+ z^kQ(z)$, for some $k\in \{1,...,n\}$, where $Q$ is a polynomial and $Q(0)\neq 0$. Substituting this equation, at $z=r\zeta$, where $r\geq 0$ and $\zeta$ is arbitrary in $\mathbb C$, in inequality (1), we arrive at
\[ 2r^k\textrm{Re}\big[\,\overline{P(0)}\zeta^kQ(r\zeta)\,\big] \ +\ r^{2k}|\zeta^kQ(r\zeta)|^2\geq 0\,, \ \forall r\geq 0\,, \ \forall \zeta \in \mathbb C,\]
and, dividing the above inequality by $r^k>0$, we find the inequality
\[\ \ \ \ \ \ \ \ \ \ \  2\textrm{Re}\big[\,\overline{P(0)}\zeta^{\,k}Q(r\zeta)\,\big] + r^{k}\big|\zeta^kQ(r\zeta)\big|^{\,2} \geq 0\,, \ \forall\, r> 0\,, \forall \,\zeta \in \mathbb C\,,\ \ \ \ \ \ \ \ \ \ \ \ \ \ \ \ \ \ \ \ \]
whose left side is a continuous function of $r$, $r\in [0,+\infty)$. 

Thus, taking the limit as $r\to 0$ we find, 
\[(2)\ \ \ \ \ \ \ \ \ \ \ \ \ \ \ \ \  \ \ \ \ \ \ \ \ \ \ \ \ \ \ \ \ \ \ 2\textrm{Re}\big[\,\overline{P(0)}Q(0)\zeta^{\,k}\,\big]\geq 0\,,\  \forall \,\zeta \in  \mathbb C\ . \ \ \ \ \ \ \ \ \ \ \ \ \ \ \ \ \ \ \ \ \ \ \ \ \ \ \ \ \ \ \ \ \ \ \ \ \ \ \ \ \ \ \]

Let $\alpha=\overline{P(0)}Q(0) = a +ib$, where $a,b \in \mathbb R$. If $k$ is odd then, substituting $\zeta=\pm 1$ and $\zeta = \pm i$ in (2), we conclude that $a=0$ and $b=0$. Hence $\alpha =0$ and then, $P(0)=0$. Thus, the case $k$ odd is proved. Next, let us suppose $k$ even. Taking $\zeta=1$ in (2), we conclude that $a\geq 0$. Picking $\zeta$ as in the lemma, let us write $\zeta^k = x+iy$, with $x<0$ and $y>0$. Substituting $\zeta^k$ and $\overline{\zeta}^k=\overline{\zeta^k}$ in (2) it follows that Re$[\alpha(x \pm iy)]=ax \mp by \geq 0$. Hence $ax\geq 0$ and (since $x<0$) we conclude that $a\leq 0$. So, $a=0$. Therefore, we get that $\mp by\geq 0$. Hence, since $y\neq 0$, we conclude that $b=0$. Hence $\alpha =0$ and then, $P(0)=0$. Thus, the case $k$ even is proved. The theorem is proved.

\

\hspace{-0,5 cm}{\bf Remarks} 

\begin{itemize}
\item[(1)] The almost algebraic ``Gauss' Second Proof'' (see [6]) of the FTA uses only that ``every real polynomial of odd degree has a real zero'' and the existence of a positive square root of every positive real number. Nevertheless, this proof by Gauss is not elementary.
\newpage
\item[(2)] It is possible to rewrite a small part of the given proof of the FTA so that the polynomial continuity is used only to guarantee the existence of $z_0$, a point of global minimum of $|P|$. In fact, to avoid extra use of  polynomial continuity, let us keep the notation of the proof and indicate $Q(z)= Q(0) + zR(z)$, with $R$ a polynomial. Then, 
substituting this expression for $Q(z)$ only in the first parcel in the left side of the inequality $ 2r^k\textrm{Re}\big[\,\overline{P(0)}\zeta^kQ(r\zeta)\,\big] \ +\ r^{2k}|\zeta^kQ(r\zeta)|^2\geq 0\,, \ \forall r\geq 0\,, \ \forall \zeta \in \mathbb C,$ that appeared just above inequality (2), we get the inequality
\[ 2\textrm{Re}\big[\,\overline{P(0)}\zeta^{\,k}Q(0)\,\big]  +  2r\textrm{Re}\big[\,\overline{P(0)}\zeta^{k+1}R(r\zeta)\,\big] + r^{k}\big|\zeta^kQ(r\zeta)\big|^{\,2} \geq 0\,,\]
$  \forall\, r> 0\,,\forall \,\zeta \in \mathbb C$. Fixing $\zeta$ arbitrary in  $\mathbb C$, it is clear that there exists $M=M(\zeta)>0$ such that $\max\big(\,|P(0)\zeta^{k+1}R(r\zeta)|, |\zeta^kQ(r\zeta)|^2\,\big)\leq M$,  $\forall r\in (0,1)$. Hence,
\[  -2\textrm{Re}\big[\,\overline{P(0)}\zeta^{\,k}Q(0)\,\big]  \leq  2rM + r^kM \leq 3rM \,, \ \forall\, r\in (0,1) .\] 
So, we conclude that $-2\textrm{Re}\big[\,\overline{P(0)}\zeta^{\,k}Q(0)\,\big] \leq 0$, with $\zeta$ arbitrary in $\mathbb C$. Now, obviously, the proof continues as in the proof of the theorem.

\item[(3)] It is worth to point out that this
proof of the FTA easily implies an independent proof of the existence of a unique positive nth root, $n\geq 3$, of any number $a\geq 0$. In fact, given $z\in \mathbb C$ such that $z^n=a$, we have that $|z|^n=a$, with $|z|\geq 0$. The uniqueness of such nth root is very trivial.

\end{itemize}

\hspace{-0,5 cm}{\bf Acknowledgments}

I thank Professors J. V. Ralston and Paulo A. Martin for very valuable comments and suggestions, J. Aragona for reference [5] and R. B. Burckel for references [7] and [13].

\

\hspace{-0,5 cm}{\bf References}

\hspace{-0,5 cm}[1] Argand, J. R., {\sl{Imaginary Quantities: Their Geometrical Interpretation}}, University Michigan Library, 2009.

\hspace{-0,5 cm}[2] Argand, J. R., {\sl{Essay Sur Une Manière de Représenter Les Quantités Imaginaires Dans Les Contructions Géométriques}}, Nabu Press, 2010.

\hspace{-0,5 cm}[3] Burckel, R. B., `` A Classical Proof of The Fundamental Theorem of Algebra Dissected'', {\sl {Mathematical Newsletter of the Ramanujan Mathematical Society 7}}, no. 2 (2007), 37-39. 

\hspace{-0,5 cm}[4] Estermann, T., ``On The Fundamental Theorem of Algebra'', {\sl{J. London Mathematical Society}} 31 (1956), 238-240.

\hspace{- 0,5 cm}[5] Fefferman, C., ``An Easy Proof of the Fundamental Theorem of Algebra'', {\sl{American Mathematical Monthly}} 74 (1967), 854-855.

\hspace{- 0,5 cm}[6] Gauss, C. F., {\sl Werke}, Volume 3, 33-56 (in latin; English translation available at
http:www.cs.man.
ac.uk/$\sim$pt/misc/gauss-web.html).

\hspace{-0,5 cm}[7] Kochol, M., ``An Elementary Proof of The Fundamental Theorem of Algebra'', International Journal of Mathematical Education in Science and Technology, 30 (1999), 614-615.

\hspace{- 0,5 cm}[8] Littlewood, J. E., ``Mathematical notes (14): Every Polynomial has a Root'', {\sl{J. London Mathematical Society}} 16 (1941), 95-98.

\hspace{-0,5 cm}[9] Oliveira, O. R. B., ``The Fundamental Theorem of Algebra: An Elementary and Direct Proof'', {\sl{ The Mathematical Intelligencer}} 33, No. 2, (2011), 1-2.

\hspace{- 0,5 cm}[10] Redheffer, R. M., ``What! Another Note Just on the Fundamental Theorem of Algebra?'', {\sl{American Mathematical Monthly}} 71 (1964), 180-185.  

\hspace{-0,5 cm}[11] Remmert,R., ``The Fundamental Theorem of Algebra''. In H.-D. Ebbinghaus, et al., {\sl{Numbers}}, Graduate Texts in Mathematics, no. 123, Springer-Verlag, New York, 1991. Chapters 3 and 4.

\hspace{-0,5 cm}[12] Searcóid, M. O., {\sl{Elements of Abstract Analysis}}, Springer-Verlag, London, 2003.

\hspace{-0,5 cm}[13] Vaggione, D., ``On The Fundamental Theorem of Algebra''. {\sl{Colloquium Mathematicum}} 73 No. 2 (1997), 193-194. 

\vspace{0,3 cm}

\hspace{- 0,5 cm}Departamento de Matemática 

\hspace{- 0,5 cm}Universidade de São Paulo

\hspace{- 0,5 cm}Rua do Matão 1010, CEP 05508-090 

\hspace{- 0,5 cm}São Paulo, SP

\hspace{- 0,5 cm}Brasil

\hspace{- 0,5 cm}e-mail: oliveira@ime.usp.br

\newpage

\end{document}